\documentclass{amsart}

\usepackage{amssymb}
\usepackage[all]{xy}
\usepackage{hyperref}
\newtheorem{thm}{Theorem}%[section]

 \newtheorem{claim}[thm]{Claim}  
   
\theoremstyle{definition}

\newtheorem{say}[thm]{}

\newtheorem*{ack}{Acknowledgments}      % \renewcommand{\theack}{} 

\newtheorem{defn-thm}[thm]{Definition--Theorem}  %!!!!!!!!!!!!!!!!!!!!!!!!
\newtheorem{defn-lem}[thm]{Definition--Lemma}  %!!!!!!!!!!!!!!!!!!!!!!!!
\theoremstyle{remark}
\renewcommand{\c}[0]{{\mathbb C}}  

\renewcommand{\o}[0]{{\mathcal O}} 
\newcommand{\z}[0]{{\mathbb Z}}

\renewcommand{\a}[0]{{\mathbb A}}

\newcommand{\p}[0]{{\mathbb P}}
\newcommand{\f}[0]{{\mathbb F}}

\newcommand{\map}[0]{\dasharrow}
\newcommand{\qtq}[1]{\quad\mbox{#1}\quad}

\newcommand{\pic}[0]{\operatorname{Pic}}

\newcommand{\rank}[0]{\operatorname{rank}}

\newcommand{\sing}[0]{\operatorname{Sing}}

\newcommand{\grass}[0]{\operatorname{Grass}}

\newcommand{\cl}[0]{\operatorname{Cl}}

\newcommand{\tsum}[0]{\textstyle{\sum}}

\def\loccoh#1.#2.#3.#4.{H^{#1}_{#2}(#3,#4)}

\DeclareMathAlphabet{\mathchanc}{OT1}{pzc}%
                                {m}{it}

\usepackage[all]{xy}\xyoption{dvips}

\begin{document}
\bibliographystyle{amsalpha}

\hfill\today

\title[Comments on Ottem--Rennemo]{Comments on: \\
  FANO VARIETIES WITH TORSION IN  $H^3(X, {\mathbb Z})$ \\
by J.~C.~OTTEM and J.~V.~RENNEMO }
  \author{J\'anos Koll\'ar}

        \begin{abstract}  Ottem and Rennemo constructed a Fano 4-fold with torsion in $H_2(X, \z)$. We simplify the computation of  $H_2(X, \z)$ and also
      exhibit 2 lines $L', L''\subset X$ such that  $L'- L''$ generates the  torsion.
  \end{abstract}

\maketitle

\cite{ottem2023fano} constructs a smooth Fano 4-fold
$X$ such that $H_2(X, \z)\cong \z+\z/2$.
%% $H^3(X, \z)\cong \z/2$ and $H^2(X, \z)\cong \z$.
Below is a shorter computation of $H_2(X, \z)$, see Claims~\ref{c1} and~\ref{c2-1}.

We also add two new results. In  Claim~\ref{c2X} we exhibit 2 lines $L', L''\subset X$ such that  $L'- L''$ generates the  torsion in  $H_2(X, \z)\cong \z+\z/2$.

Then 
in Paragraph~\ref{s15} we show that $X$ is birational to a double cover of $\p^4$ ramified along a degree 18 hypersurface $R$, which is obtained as the 5-secants of a degree 15, smooth,  determinantal surface
$
S=\bigl(\rank N\leq 4\bigr)\subset \p^4,
$
where $N$ is a  $6\times 5$  matrix
whose entries are linear forms.
Although $S$ is smooth, it  is not a general determinantal surface, since
the latter have only  1-parameter families of 5-secants.

The higher dimensional examples constructed in \cite{ottem2023fano}
can  also be treated with minor changes.

We refer to \cite[Chaps.VIII-IX]{Room-Det} for  symmetric determinantal varieties and to \cite{zbMATH02689354} for the classification of lines on them. 
\medskip

\begin{say}[Basic set-up]\label{s1}  Following \cite{ottem2023fano},
let $Z_i\subset \p^{14}$ be the space of rank $\leq i$ quadrics in $\p^4$.
Our main interest is 
$Z=Z_4\subset \p^{14}$, the space of rank $\leq 4$ quadrics in $\p^4$.
It is a quintic hypersurface.

The universal deformation of a rank 3 quadric is given by
$$
x_0^2-x_1^2-x_2^2+z_1x_3^2+z_2x_3x_4+z_3x_4^2=0.
\eqno{(\ref{s1}.1)}
$$
This has rank $\leq 4$ iff $z_2^2-4z_1z_3=0$. 
So $Z$   is singular along
$Z_3$,   with transversal singularity type $A_1$.  

Let $Z^\circ:= Z_4\setminus Z_2$ be the space  of quadrics of rank $3$ or $4$.
Note that $Z_2\subset Z_4$ has codimension 5.

Let $U$ be the  space of pairs  $(L^2\subset Q)$ where $L^2\subset\p^4$ is a 2-plane and $Q\subset\p^4$ a quadric. Note that this forces $\rank Q\leq 4$, and
$U$ is a $\p^8$-bundle over $G:=\grass(2,4)$.
We have projections
$$
G\stackrel{p_1}{\longleftarrow} U \stackrel{p_2}{\longrightarrow} Z,
$$
and $\pic(U)\cong \z^2$ is generated by $p_1^*\o_G(1)$ and $p_2^*\o_{Z}(1)$. 

Let $U_2\subset U$ be the subset where the quadric has rank $\leq 2$.
These quadrics split into 2 hyperplanes, one of which must contain $L^2$. 
So $U_2$  is a $\p^1\times \p^4$-bundle over $\grass(2,4)$.
Set $U^\circ:=U\setminus U_2$.

Since $U_2\subset U$ has codimension 3, 
$H^i(U^\circ, \z)=H^i(U, \z)$ for $i\leq 4=2(3-1)$.
In particular, 
$$
H^0(U^\circ, \z)\cong \z, \ 
H^1(U^\circ, \z)\cong 0, \ 
H^2(U^\circ, \z)\cong \z^2, \ 
H^3(U^\circ, \z)\cong 0.
$$

Now look at the morphism  $p_2: U\to Z$ with  Stein factorization
$$
U\stackrel{\pi}{\longrightarrow} W \stackrel{\tau}{\longrightarrow}Z.
$$
In the coordinates (\ref{s1}.1),  $Z\cong (z_2^2-4z_1z_3=0)\times \a^{11}$.
Set $w_1:=\sqrt{z_1}$ and $w_3:=\sqrt{z_3}$.  Since $w_1w_3=z_2/2$, adjoining both $w_1, w_2$ is a degree 2 covering and 
$$
z_1x_3^2+z_2x_3x_4+z_3x_4^2=(w_1x_3+w_3x_4)^2.
$$
Thus, locally analytically over the points of
$Z_3\setminus Z_2$,  we have 
$W\cong  \a^2_{w_1, w_2}\times \a^{11}$, and the family of quadrics becomes
$$
(x_0^2-x_1^2-x_2^2+(w_1x_3+w_2x_4)^2=0)\times \a^{11}.
\eqno{(\ref{s1}.2)}
$$
Thus $U$ is the family of 2-planes in the same family as
$$
\bigl(x_0-x_1=x_2-(w_1x_3+w_2x_4)=0\bigr).
$$
Each of these has a unique intersection point with
$(x_3=x_4=0)$. Thus   $U$  is locally analytically
isomorphic to the trivial family
$$
W\times (x_0^2-x_1^2-x_2^2=0)\subset W\times \p^2.
$$
Therefore, 
restricting to $U^\circ $ we get 
$$
U^\circ\stackrel{\pi^\circ}{\longrightarrow} W^\circ \stackrel{\tau^\circ}{\longrightarrow}Z^\circ,
$$
and   $\pi^\circ: U^\circ\to W^\circ$ is a smooth morphism with conics as fibers.
\end{say}

\begin{claim}\label{c1} $H^3(W^\circ, \z)\cong \z/2$ and $\pi$ has  no rational sections.
\end{claim}

Proof. Let $C\subset U^\circ$ be a fiber of $\pi^\circ$ over a rank 3 conic.
We can think of $C$ as (cones over) lines on a  quadric cone, or after further degeneration, as  (cones over) 2 pencils of lines on 2 planes.
So $C$ is 2-times the class of (cones over) a  pencil of lines in a plane.
Thus the image of $H^2\bigl(  U^\circ, \z\bigr)\to H^2(C, \z)\cong \z$ is twice the generator.
(Note that this splitting of $C$ into 2 components happens in the fibers over $Z_2$, thus outside  $Z^\circ$.)

In the Leray spectral sequence for $\pi^\circ$, the only interesting map is on the $E_3$ page:
$$
\z \cong H^0\bigl(  W^\circ, R^2\pi^\circ_*\z\bigr) \to H^3\bigl(  W^\circ, \z\bigr).
$$
As we noted, the kernel is $2\z$ and $ H^3\bigl(  U^\circ, \z\bigr)=0$.
Thus $ H^3\bigl(  W^\circ, \z\bigr)\cong \z/2$  and $H^2\bigl(  W^\circ, \z\bigr)\cong \z$.  \qed
\medskip

\begin{say}[Construction of $X$] \label{s3}
Since $Z\subset \p^{14}$ is a quintic hypersurface,
$K_Z\sim -10H$. Since  $\tau: W\to Z$ is \'etale over $Z\setminus Z_3$,
$K_W\sim -10 \tau^*H$.
As we noted, $Z_2\subset Z$ has codimension $5$, so its preimage 
$W_2\subset W$ has codimension $5$.

As in \cite{ottem2023fano},
let $X_Z\subset Z$ be the complete intersection of $9$ general hyperplanes
and  $X:=X_W\subset W$ its preimage.
Then $X\subset W^\circ$ and $K_X\sim -\tau^*H$.
\end{say}

\begin{claim}\label{c2-1} $H_2(X, \z)\cong \z+\z/2$.
\end{claim}

Proof. The isomorphism $H_2(X, \z)\cong H_2(W^\circ, \z)$ follows from
the Lefschetz hyperplane theorem. A slight problem is that $W$ is singular and
$W^\circ$ is not compact. There are several ways to take care of these.

One way is to use  the Lefschetz theorem in intersection homology.
A more elementary argument is given in \cite[4.9]{ottem2023fano}. \qed

\medskip

Note that by \cite{MR2306166}, $H_2(X, \z)$ is generated by algebraic curves.
Next we write down a difference of 2 smooth, degree 1 rational curves that generates
the $\z/2$-summand of $H_2(X, \z)$.
\medskip

\begin{claim}\label{c2}  Let $L\subset Z\setminus Z_3\subset \p^{14}$ be a line.
  Its preimage in $W$ is a pair of lines $L'\cup L''$ such that
  \begin{enumerate}
\item  $L'$ and $L''$ are numerically equivalent, 
  \item  $U\times_WL'$ is the ruled surface $\f_1\cong B_p\p^2$,
  \item  $U\times_WL''$ is the ruled surface $\f_0\cong \p^1\times \p^1$,  and
    \item  $L'-L''$  is a  generator of
      the $\z/2$-summand  of $H_2(W^\circ, \z)\cong \z+\z/2$.
      \end{enumerate}
\end{claim}

By Paragraph~\ref{lines.XZ.say}, 
$X_Z$   contains a
2-parameter family of lines, and  Claim~\ref{c2} applies to them.
Thus we obtain the following.

\begin{claim}\label{c2X}  Let $L\subset X_Z\setminus Z_3$ be a line.
Its preimage in $X$ is a pair of lines $L'\cup L''$, and
$L'-L''$  is a  generator of
the $\z/2$-summand  of $H_2(X, \z)\cong \z+\z/2$. \qed
\end{claim}

\begin{say}[Beginning of the proof of Claim~\ref{c2}]
By Claim~\ref{c2-1}  $H_2(X, \z)/(\mbox{torsion})$ is  generated by $K_X\sim -\tau^*H$,  
 so $L'$ and $L''$ are numerically equivalent.

  By Paragraph~\ref{lines.say},   in suitable coordinates  we can write
$L$ as a family of quadrics
$$
Q(\lambda{:}\mu):=\bigl(x_0(\lambda x_2-\mu x_3)=x_1(\mu x_4-\lambda x_3)\bigr).
$$
All of these contain the 2-plane  $(x_0=x_1=0)$,  defining  a section
$s_0:L\to U$.

The preimage of $L$ in $W$ is a disjoint union of 2 lines
$L'\cup L''$.  We choose $L'$ to be
$\pi\circ s_0:L\to U\to W$.

For any nonzero linear form  $\ell=a\lambda+b\mu $, a section  $C'(\ell)$ of
$\pi: U\to W$ over $L'$ is given by
$$
\ell x_0=\mu x_4-\lambda x_3 \qtq{and} \lambda x_2-\mu x_3=\ell x_1.
$$
For $\ell_1\neq \ell_2$ the 2 sections  $C'(\ell_1), C'(\ell_2)$ meet at the point where $\ell_1=\ell_2$.
Thus  $U\times_WL'$ is the ruled surface  $\f_1$.

In the other family of 2-planes, we have sections  $C''(\ell)$ given by
$$
cx_0=x_1 \qtq{and} \mu x_4-\lambda x_3=c(\lambda x_2-\mu x_3).
$$
These are disjoint for $c_1\neq c_2$. Thus $U\times_WL''$ is the the trivial $\p^1$-bundle.

These show  Claim~\ref{c2}.2--3. \qed
\end{say}

Claim~\ref{c2}.4 is a formal consequence of (\ref{c2}.1--3).
To see this, 
 we need to  discuss how to detect 2-torsion in $H_2$ using
$\p^1$-bundles.  (Similarly, $n$-torsion can be detected using
$\p^{n-1}$-bundles.)

\begin{say}[Comments on $\p^1$-bundles]\label{com1}
Let $X$ be a  normal,  proper variety and 
$\pi:Y\to X$ a $\p^1$-bundle  (\'etale locally trivial).
For a smooth curve $C\to X$, let    $C_Y\to Y_C:=C\times_XY$ be a lifting. Set
$$
\sigma_Y(C):=(C_Y\cdot K_{Y_C/C}) \mod 2.
\eqno{(\ref{com1}.1)}
$$
This is well defined as a function on $A_1(X)$, the group of 1-cycles modulo algebraic equivalence.

If $X$ is smooth and $Y\to X$ has  a rational section  $S\subset Y$, then
$ K_{Y/X}\sim -2S+\pi^*D$ for some Cartier divisor $D$ on $X$. In this case
$$
\sigma_Y(C)\equiv (C\cdot D) \mod 2.
$$
Conversely, assume that  we are over $\c$,  $H_2(X,\z)$ is generated by algebraic curves, and $ \sigma_Y(C)\equiv (C\cdot D) \mod 2$  for every $C$. Then 
$K_{Y/X}-\pi^*D$ is divisible by 2, giving a rational section.
In any case, we get the following.
\end{say}

\begin{claim}\label{c3} If there is a numerically trivial 1-cycle $C$ such that
$\sigma_Y(C)\equiv 1 \mod 2$, then
$[C]$ is a nontrivial torsion element in  $H_2(X,\z)$,
and $Y\to X$ has no rational sections. \qed
\end{claim}

\begin{say}[End of the proof of Claim~\ref{c2}]
  Since   $U\times_WL'\cong \f_1$, (\ref{com1}.1)   shows that
   $\sigma_U(L')\equiv 1 \mod 2$. Similarly, 
  $U\times_WL''\cong \f_0$ implies that 
  $\sigma_U(L'')\equiv 0 \mod 2$.
  Thus $C:=L'-L''$ is numerically trivial and
   $\sigma_U(C)\equiv 1 \mod 2$.
We can now apply Claim~\ref{c3}. \qed
\medskip

In both cases we could have used the isomorphism
$$
\omega_{U/W}\cong p_1^*\o_{G}(-1)\otimes p_2^*\o_{Z}(1)
$$
to compute $\sigma_U(L')$ and $\sigma_U(L'')$. \qed
\end{say}

\begin{say}[Lines on $Z$]\label{lines.say} By  \cite{zbMATH02689354}, the lines on $Z\setminus Z_2$ form 3 families of dimension 20 each. These are the following.

\medskip
(\ref{lines.say}.1) $\langle Q_1, Q_2\rangle$   where the $Q_i$ contain a common 2-plane. The general such line $L$ is disjoint from $Z_3$, and its preimage in $W$ is a pair of disjoint lines $L'\cup L''$. After coordinate change, these can be written as
$$
x_0(\lambda x_2-\mu x_3)=x_1(\mu x_4-\lambda x_3).
$$

\medskip
(\ref{lines.say}.2)  $\langle Q_1, Q_2\rangle$   where the $Q_i$ have a common singular point. After coordinate change, these can be written as
$$
\lambda q_1(x_1, \dots, x_4)+\mu q_2(x_1, \dots, x_4)=0,
$$
where the $q_i$ are quadratic forms.
The general such line $L$ intersects $Z_3$ at 4 points,  and its preimage in $W$ is a smooth, elliptic curve of degree 2.

\medskip
(\ref{lines.say}.3)  $\langle Q_1, Q_2\rangle$   where the $Q_i$ have rank 2 and $\sing Q_i$ is tangent to $Q_{3-i}$.  The general such line $L$ intersects $Z_3$ at 2 points,  and its preimage in $W$ is a smooth, rational curve of degree 2. After coordinate change, these can be written as
$$
\lambda(x_0^2+x_1x_2)+\mu(x_2x_3+x_4^2)=0.
$$
\end{say}

\begin{say}[Lines on $X_Z$]\label{lines.XZ.say}
The space of lines in $Z\setminus Z_3$ has dimension 20, and with each hyperplane section the dimension drops by 2. So the lines on $X_Z$ form 3 families of dimension 2 each. Only (\ref{lines.say}.1) contains lines that are disjoint from $\sing X_Z$. 

Since there are no lines on $Z_3\setminus Z_2$, 
the only common liness to any 2 of these families are the finitely many double tangents of $\sing X_Z$. 
\end{say}

\begin{say}[Another representation of $X_Z$]\label{s10}
  Let $P_5$ be a general 5-dimensional linear system of quadrics on $P_4:=\p^4$.
For $i=4,5$, we have the projections $\pi_i:P_4\times P_5\to P_i$. For brevity let us write  $(a,b):=a\pi_4^*H_4+ b \pi_5^*H_5$, where
 $H_i$ is the hyperplane class on $P_i$.
  Set
  $$
  Y:=\bigl\{(p, Q): p\in P_4, Q\in P_5, p\in \sing Q\bigr\}\subset P_4\times P_5.
  $$
  The condition $ p\in \sing Q$ is equivalent to the partial derivatives of the equation of $Q$ vanishing at $p$. Thus $Y\subset P_4\times P_5$ is the complete intersection of 5 divisors of bidegree $(1,1)$. Write these as
  $$
  \tsum_{i=0}^4\tsum_{j=0}^5 a^\ell_{ij}y_ix_j\qtq{for} \ell=1,\dots, 5.
  \eqno{(\ref{s10}.1)}
  $$
  Over $P_5$, (\ref{s10}.1) is equivalent to a $5\times 5$ symmetric matrix
  whose entries are the linear forms  $m_i^\ell=\tsum_{j=0}^5 a^\ell_{ij}x_j$.
  The condition $\det (m_i^\ell)=0$ defines  $X_Z$ as in Paragraph~\ref{s3}.

  Over $P_4$, (\ref{s10}.1) is equivalent to a $6\times 5$  matrix
  whose entries are the linear forms  $n_j^\ell=\tsum_{i=0}^4 a^\ell_{ij}y_i$.
  Note that $\pi_4:Y\to P_4$ is birational. Its inverse is  the blow-up
  of a surface
  $$
  S\subset P_4,\qtq{defined by }\rank(n_j^\ell)\leq 4.
  $$
  Let $E_4\subset Y$ denote the exceptional divisor.
   $Y$ defines a rational map $P_4\map X_Z\subset P_5$, which is given by the
  $5\times 5$ subdeterminants of $(n_j^\ell)$. Thus
  $E_4\sim (5,-1)|_Y$.

  The inverse rational map  $X_Z\map P_4$ is a bit harder to see.
  It is given by a linear system of divisors as follows.
  Let $H\in |H_4|$ be a hyperplane and set
  $$
  D_H:=\bigl\{Q\in X_Z: H\cap \sing Q\neq\emptyset\bigr\}\subset X_Z.
  $$
  Note that the condition $H\cap \sing Q$ is equivalent to
  $Q|_H$ being singular. (Here we need that $Q$ itself is singular.)
  This gives us the equation  $\det (Q|_H)=0$ for $D_H$. It has degree 4.

  We claim that the intersection of $\left(\det (Q|_H)=0\right)$ with $X_Z$ has multiplicity 2.
  To see this, choose coordinates such that  $H=(x_0=0)$ and
  $Q=(x_0^2+x_2^2+x_3^2+x_4^2=0)$. For its deformations we can make linear coordinate changes to the $x_2, x_3, x_4$, but $x_0$  can only be multiplied by a constant. Thus we get a miniversal deformation family
  $$
  \bigl(x_0^2+t_1x_0x_1+t_2x_1^2+x_2^2+x_3^2+x_4^2=0\bigr)\subset 
P_5\times \a^2_{t_1, t_2}.
  $$
  For a given $t_1, t_2$, the quadric has rank $4$ iff $t_1^2-4t_2=0$, and the singular point is on $(x_0=0)$ iff $t_2=0$.  Their intersection is the length 2 scheme $(t_1^2=0)$.

  Thus the $D_H\subset  X_Z$ have degree $10=\frac12 (4\cdot 5)$ and $2D_H\sim 4H_5|_Y$.
  In particular, the divisor class $D_H-2H_5|_Y$ is 2-torsion in the class group $\cl(X_Z)$.
  The corresponding double cover is our  $X$, constructed  in Paragraph~\ref{s3}.

  Let $E_5\subset Y$ denote the exceptional divisor of $\pi_5:Y\to X_Z$.
  The previous computations suggest that it should be linearly equivalent to $(-1,2)$. However, $X_Z$ has multiplicity 2 along the base locus of
  $|D_H|$, so the correct bidegree is $(-2,4)$.

On $P_4$, the 3 families of  lines  (\ref{lines.say}.1--3)  correspond to

(\ref{s10}.1) conics that are 9-secants of $S$,

(\ref{s10}.2) fibers of $E_4\to S$, and 

(\ref{s10}.3) lines that are 4-secants of $S$.

\end{say}

\begin{say}[$X$ as a double $\p^4$]\label{s15}
  By the previous description, $X$ is birational to a double cover of $\p^4$ ramified along the hypersurface  $R:=\pi_4(E_5)\subset P_4$.

  The degree of $R$ is given by  $(1,0)^3[E_5](1,1)^5$, which works out to be 18.
  The degree of the surface $S$ is $(1,0)^2[E_4](0,1)(1,1)^5=15$.
  Note that $S$ is a $6\times 5$ determinantal surface. However, it is not general since we have a symmetry condition on the  $P_5$ side, so  results about general  determinantal surfaces do not apply to $S$.

  The multiplicity of $Y$ along $S$ is 4. This follows from the computation
  $$
  (1,0)^2[E_4][E_5](1,1)^5=60 = 4\cdot \deg S.
  $$
Thus $R$ is in the 4th symbolic power of the homogeneous ideal of $S$, but not in its 4th power. For general  determinantal surfaces these  are equal by \cite{MR0558865}.

  Another interesting property of $S$ is that the fibers of $E_5\to \sing X_Z$
  give 5-secants of $S$. Thus $S$ has a 2-parameter family of 5-secants.
Note that  by (\ref{s10}.3), the family of  4-secants   has  dimension 2 as well.

  Most surfaces in $\p^4$, including  general $6\times 5$ determinantal surfaces,  have only  1-parameter families of 5-secants.

  \end{say}

\begin{ack}  I thank  J.~Ottem, C.~Raicu and B.~Totaro   for many   useful comments.
Partial  financial support    was provided  by  the NSF under grant number
DMS-1901855.
\end{ack}

%%\bibliography{../refs-main/refs}

\begin{thebibliography}{dCEP80}

\bibitem[dCEP80]{MR0558865}
C.~de~Concini, D.~Eisenbud, and C.~Procesi, \emph{Young diagrams and
  determinantal varieties}, Invent. Math. \textbf{56} (1980), no.~2, 129--165.
  \MR{558865}

\bibitem[Kro1890]{zbMATH02689354}
L.~Kronecker, \emph{Algebraische {Reduction} der {Scharen} bilinearer
  {Formen}}, Berl. Ber. \textbf{1890} (1890), 1225--1237.

\bibitem[OR23]{ottem2023fano}
J.~C. Ottem and J.~V. Rennemo, \emph{Fano varieties with torsion in the third
  cohomology group}, 2023, \url{https://arxiv.org/abs/2309.10793}.

\bibitem[Roo38]{Room-Det}
T.~G. Room, \emph{The geometry of determinantal loci}, Cambridge University
  Press, London, 1938.

\bibitem[Voi06]{MR2306166}
Claire Voisin, \emph{On integral {H}odge classes on uniruled or {C}alabi-{Y}au
  threefolds}, Moduli spaces and arithmetic geometry, Adv. Stud. Pure Math.,
  vol.~45, Math. Soc. Japan, Tokyo, 2006, pp.~43--73. \MR{2306166}

\end{thebibliography}

\def\cprime{$'$} \def\cprime{$'$} \def\cprime{$'$} \def\cprime{$'$}
  \def\cprime{$'$} \def\dbar{\leavevmode\hbox to 0pt{\hskip.2ex
  \accent"16\hss}d} \def\cprime{$'$} \def\cprime{$'$}
  \def\polhk#1{\setbox0=\hbox{#1}{\ooalign{\hidewidth
  \lower1.5ex\hbox{`}\hidewidth\crcr\unhbox0}}} \def\cprime{$'$}
  \def\cprime{$'$} \def\cprime{$'$} \def\cprime{$'$}
  \def\polhk#1{\setbox0=\hbox{#1}{\ooalign{\hidewidth
  \lower1.5ex\hbox{`}\hidewidth\crcr\unhbox0}}} \def\cdprime{$''$}
  \def\cprime{$'$} \def\cprime{$'$} \def\cprime{$'$} \def\cprime{$'$}
\providecommand{\bysame}{\leavevmode\hbox to3em{\hrulefill}\thinspace}
\providecommand{\MR}{\relax\ifhmode\unskip\space\fi MR }
\providecommand{\MRhref}[2]{%
  \href{http://www.ams.org/mathscinet-getitem?mr=#1}{#2}
}
\providecommand{\href}[2]{#2}

\bigskip

  Princeton University, Princeton NJ 08544-1000, \

  \email{kollar@math.princeton.edu}

\end{document}